\documentclass[a4paper,11pt]{article}
\usepackage[latin1]{inputenc}
\usepackage[T1]{fontenc}

\usepackage{amsfonts}
\usepackage{amsmath,amssymb}
\usepackage{pstricks}
\usepackage{amsthm}
\usepackage[all]{xy}
\usepackage[pdftex]{graphicx}
\usepackage{graphicx}   
\newtheorem{theo}{Theorem}[section]

\newtheorem{prop}[theo]{Proposition}

\newtheorem{cor}[theo]{Corollary}
\newtheorem{lem}[theo]{Lemma}
\newtheorem{rem}[theo]{Remark}
\newtheorem{ex}{Example}
\newenvironment{dem}{\noindent
  \textit{{Proof}} : }
  {\hfill\qedsymbol\newline}

\title{Jet schemes of complex plane branches and equisingularity}
\author{Hussein MOURTADA}

\begin {document}

\maketitle  

\begin{abstract}
For $m \in \mathbb{N}$, we give formulas for the number $N(m)$ of irreducible components of the m-th Jet Scheme of a complex branch $C$
and for their codimensions, in terms of $m$ and the generators of
the semigroup of $C$. This structure of the Jet Schemes determines and is determined by the topological type of $C$.
\end{abstract}

\section{Introduction}

$~~~$Let $k$ be an algebraically closed field.  The space of arcs $X_{\infty}$ of an algebraic $k-$variety $X$ is a non-noetherian
scheme in general. It has been introduced by Nash in \cite{N}. Nash has initiated its study by looking at its image by the truncation maps $X_{\infty} \longrightarrow X_m$ in the jet schemes of $X.$The $m^{th}-$jet scheme $X_m$ of $X$ is a $k-$ scheme
of finite type which parmametizes morphisms $Spec~\frac{k[t]}{t^{m+1}} \longrightarrow X.$ From now on we assume $char~k=0.$
In \cite{N}, Nash has derived from the existence of a resolution of singularities of $X$, that the number of irreducible components
of the Zariski closure of the set of the $m-$truncations of arcs on $X$ that send $0$ into the singular locus of $X$ is constant for $m$ large enough.
Besides a theorem of Kolchin asserts that if $X$ is irreducible, then $X_{\infty}$ is also irreducible. More recently ,
the jet schemes have attracted attention from various viewpoints. In \cite{M},Mustata has characterized the locally complete intersection varieties having irreducible $X_m$ for $m\geq 0.$In \cite{ELM} , a formula comparing the codimensions of $Y_m$ in $X_m$ with the log canonical threshold of a pair $(X,Y)$ is given.In this work, we consider a curve $C$ in the complex plane 
$\mathbb{C}^2$ with a singularity at $0$ at which it is analytically irreducible (i.e. the formal neighborhood$(C,0)$ of $C$ at $0$ is a branch). We determine the irreducible components of the space $C_m^0:=\pi_m^{-1}(0)$ where $\pi_m:C_m\longrightarrow C$ is the canonical projection, and we show that their number is not bounded as $m$ grows. More precisely, let $x$ be a transversal parameter in the local ring 
$O_{\mathbb{C}^2,0},$ i.e. the line $x=0$ is transversal to $C$ at $0$ and following \cite{ELM},for $e \in \mathbb{N}$ let  
$$Cont^e(x)_{m}(resp.Cont^{>e}(x)_{m}):=\{ \gamma \in C_m \mid ord_tx \circ \gamma=e(resp.>e)\}.$$
 
Let $\Gamma(C)=< \overline{\beta}_0, \cdots ,\overline{\beta}_g>$ be the semigroup of the branch $(C,0)$ and let
$e_i=gcd(\overline{\beta}_0, \cdots ,\overline{\beta}_i),$  $0\leq i\leq g.$ Recall that $\Gamma(C)$ and the topological type of $C$ near $0$ are equivalent data. We show in theorem $4.9$ that the irreducible components of $C_m^0$ are

$$C_{m \kappa I}=\overline{Cont^{\kappa\bar{\beta}_0}(x)_m},$$for   $1 \leq \kappa$ and $\kappa\bar{\beta}_0\bar{\beta_1}+e_1 \leq m,$
 $$C^j_{m\kappa v}=\overline{Cont^{\frac{\kappa \bar{\beta}_0}{e_{j-1}}}(x)_m} $$ 
for $2\leq j \leq g,1 \leq \kappa, \kappa \not \equiv 0~ mod ~\frac{e_{j-1}}{e_j}$ and $\kappa \frac{\bar{\beta}_0\bar{\beta}_1}{e_{j-1}}+e_1\leq m< \kappa\bar{\beta}_j,$\
 
$$B_m= Cont^{>n_1q}(x)_m,$$
if $qn_1\bar{\beta_1}+e_1 \leq m <(q+1)n_1\bar{\beta}_1+e_1.$\\
These irreducible components give rise to infinite and finite inverse systems represented by a tree.We recover $<\overline{\beta}_0, \cdots ,\overline{\beta}_g>$
from the tree and the multiplicity $\overline{\beta}_0$ in corollary $4.13,$ and we give formulas for the number of irreducible components of $C_m^0$ and their codimensions in terms of $m$ and $(\overline{\beta}_0, \cdots ,\overline{\beta}_g)$ in proposition $4.7$ and corollary $4.10.$ We recover the fact coming from \cite{ELM} and \cite{I} that 
$$min_m \frac{codim(C_m^0,\mathbb{C}^2_m)}{m+1}=\frac{1}{\overline{\beta}_0} +\frac{1}{\overline{\beta}_1}.$$

$~~~$ The structure of the paper is as follows:
The basics about Jet schemes and the results that we will need are presented in section 2. In section 3 we present the definitions
and the reults we will need about branches. The last section is devoted to the proof of the main result and corollaries.\\
\\ AKNOWLEDGEMENTS\\
\\
 I would like to express all my gratitude to Monique Lejeune-Jalabert, without whom this work would not exist.

\section{Jet schemes}
Let $k$ be an algebraically closed field of arbitrary characteristic.
Let $X$ be a $k$-scheme of finite type over $k$ and let $m \in \mathbb{N}$. The functor $F_m :k-Schemes \longrightarrow Sets$
which to an affine scheme defined by a $k-$algebra $A$ associates
\\ $$F_m(Spec(A))=Hom_k(Spec A[t]/(t^{m+1}),X)$$
is representable by a $k-$scheme $X_m$ \cite{V}. $X_m$ is the m-th jet scheme of $X$, and $F_m$ is  isomorphic to its functor of points.
In particular the closed points of $X_m$ are in bijection with the $k[t]/(t^{m+1})$ points of $X$.
\\For $m,p \in \mathbb{N}, m > p$, the truncation homomorphism $A[t]/(t^{m+1}) \longrightarrow A[t]/(t^{p+1})$ induces
a canonical projection $\pi_{m,p}: X_m \longrightarrow X_p.$ These morphisms clearly verify $\pi_{m,p}\circ \pi_{q,m}=\pi_{q,p}$
for $p<m<q$.
\\ Note that $X_0=X$. We denote the canonical projection $\pi_{m,0}:X_m\longrightarrow X_0$ by $\pi_{m}$.

\begin{ex}Let $X=Spec ~\frac{k[x_0, \cdots ,x_n]}{(f_1, \cdots ,f_r)}$ be an affine $k-$scheme. For a $k$-algebra $A$, to give a
$A$-point of $X_m$ is equivalent to give a $k$-algebra homomorphism
$$\varphi : \frac{k[x_0, \cdots ,x_n]}{(f1, \cdots ,fr)}\longrightarrow A[t]/(t^{m+1}).$$
The map $\varphi$ is completely determined by the image of $x_i,i=0, \cdots ,n$
$$x_i \longmapsto \varphi(x_i)=x_i^{(0)}+x_i^{(1)}t+  \cdots  + x_i^{(m)}t^m$$
such that $f_l(\phi(x_0),  \cdots  ,\phi(x_n))\in (t^{m+1})$,  $l=1, \cdots ,r.$
\\
\\ If we write $$f_l(\phi(x_0),  \cdots  ,\phi(x_n))=\sum_{j=0}^m F_l^{(j)}(\underline{x}^{(0)}, \cdots ,\underline{x}^{(j)})~t^jmod ~~ (t^{m+1})$$
where $\underline{x}^{(j)}=(x^{(j)}_0, \cdots ,x^{(j)}_n)$, then
$$X_m=Spec \frac{k[\underline{x}^{(0)}, \cdots ,\underline{x}^{(m)}]}{(F_l^{(j)})_{l=1, \cdots ,r}^{ j=0, \cdots ,m}} $$
\end{ex}
\begin{ex}
From the above example, we see that the m-th jet scheme of the affine space $\mathbb{A}_k^n$ is
isomorphic to $\mathbb{A}_k^{(m+1)n}$ and that the projection $\pi_{m,{m-1}}:(\mathbb{A}_k^n)_m\longrightarrow (\mathbb{A}_k^n)_{m-1}$
is the map that forgets the last $n$ coordinates.
\end{ex}
\begin{lem}
If $f:X\longrightarrow Y$ is an \'etale morphism, then for every $m \in \mathbb{N}$,
the following diagram
\\
\[  \xymatrix{
      X_m  \ar[d]_{\pi_m} \ar[r]^{f_m}  & Y_m \ar[d]^{\pi_m} \\
      X  \ar[r]_f                                      & Y                  }
\]
is cartesian.
\end{lem}
\begin{dem} For a $k$-algebra $A$, to give an $A$-point of $Y_m \times_Y X$ is equivalent to give a commutative diagram
\\
\[  \xymatrix{
      Spec(A) \ar[d] \ar[r]  & X \ar[d]^{f} \\
      Spec(A[t]/(t^{m+1}))  \ar[r]      & Y                  }
\]
  which is equivalent  to give a unique morphism $Spec(A[t]/(t^{(m+1)}))\longrightarrow X$
   making the two triangles commutative,since $f$ is formally \'etale.
\end{dem}
\begin{cor}If $X$ is a nonsingular $k-$variety of dimension $n$, then all projections $\pi_{m,{m-1}}: X_m\longrightarrow X_{m-1}$
are locally trivial fibrations with fiber $\mathbb{A}_k^n$. Then in particular $X_m$ is a nonsingular variety of dimension $(m+1)n$.
\end{cor}
\begin{dem}
It is sufficient to prove that for every $x \in X$ there exists an open neighborhood $U$ of $x$ such that $U_m\simeq U\times_k \mathbb{A}_k^{mn}.$ But since $X$ is nonsingular,
there exists an open neighborhood $U$ of $x$ and an \'etale morphism $g:U\longrightarrow \mathbb{A}_k^n$. Then we deduce the claim from the above lemma .
\end{dem}

Let $char(k)=0$,  $S=k[x_0, \cdots .,x_n]$ and
$S_{m}=k[\underline{x}^{(0)}, \cdots .,\underline{x}^{(m)}]$.
Let $D$ be the $k-$derivation on $S_m$ defined by $D(x_i^{(j)})=(j+1)x_i^{(j+1)}$ if $0\leq j<m$, and $D(x_i^{(m)})=0$. For $f \in S$
let $f^{(1)}:=D(f)$ and we recursively define $f^{(m)}=D(f^{(m-1)})$. 

\begin{prop} Let $X=Spec(S/(f_1, \cdots ,f_r))=Spec(R)$ and $R_m=\Gamma (X_m)$. Then
$$R_m=Spec(\frac{k[\underline{x}^{(0)}, \cdots .,\underline{x}^{(m)}]}{(f_i^{(j)})_{i=1, \cdots ,r}^{j=0, \cdots ,m}}.$$
 
\end{prop}
\begin{dem}
For a $k-$algebra $A$, to give an $A-$point of $X_m$ is equivalent to give an homomorphism
$$ \phi:k[x_0, \cdots .,x_n]\longrightarrow A[t]/(t^{m+1})$$
which can be given by $$x_i\longrightarrow \frac{x_i^{(0)}}{0!}+\frac{x_i^{(1)}}{1!}t+  \cdots  + \frac{x_i^{(m)}}{m!}t^m.$$
Then for a polynomial $f \in S$, we have
$$\phi(f)=\sum_{j=0}^m \frac{f^{(j)}(\underline{x}^{(0)}, \cdots ,\underline{x}^{(j)})}{j!}~t^j.$$
To see this, it is sufficient to remark that it is true for $f=x_i$, and that both sides of the
equality are additive and multiplicative in $f$, and the proposition follows.
\end{dem}
\begin{rem} Note that the proposition shows the linearity of the equations $F_i^{j}(\underline{x}^{(0)}, \cdots ,\underline{x}^{(j)})$ defining $X_m$
with respect to the new variables i.e $\underline{x}^{(j)}$, which is the algebraic point of view on the fibration in corollary $2.2$.
\end{rem}
\section{ Semigroup of complex branches}
The main references for this section are \cite{Z},\cite{Me},\cite{A},\cite{Sp},\cite{GP},\cite{GT},\cite{LR}.
Let $f \in \mathbb{C} [[x,y]]$ be an irreducible power series, which is $y$-regular (i.e $f(0,y)=y^{\beta_0}u(y)$
where $u$ is invertible in $\mathbb{C}[[y]]$) and such that $mult_0f=\beta_o$  and let $C$ be the analytically irreducible plane
curve(for short branch) defined by $f$ in $Spec~\mathbb{C} [[x,y]]$.
By the Newton-Puiseux theorem, the roots of $f$ are 
$$y=\sum_{i=0}^{\infty}a_iw^ix^{\frac{i}{\beta_o}}~~~~~~~(1)$$
where $w$ runs over the $\beta_0-th$-roots of unity in $\mathbb{C}$.This is equivalent to the existence of a parametrization
of $C$ of the form $$ x(t)=t^{\beta_0}$$  $$y(t)= \sum_{i\geq \beta_0}a_i t^i.$$
We recursively define  $\beta_i=min \lbrace i, a_i \neq 0,~ gcd(\beta_0, \cdots ,\beta_{i-1})$ is not a divisor of $i \rbrace$. \\
Let $e_0=\beta_0$ and $e_i=gcd(e_{i-1},\beta_i), i\geq 1$. Since the sequence of positive integers
$$e_0>e_1> \cdots >e_i> \cdots $$
is strictly decreasing, there exists $g \in \mathbb{N}$, sucht that $e_g=1$.
The sequence $(\beta_1, \cdots .,\beta_g)$ is the sequence of Puiseux exponents of $C$. We set
$$n_i:=\frac{e_{i-1}}{e_i}, m_i:=\frac{\beta_i}{e_i}, i=1, \cdots ,g$$
and by convention, we set $\beta_{g+1}=+\infty$ and $n_{g+1}=1$.
\\
\\On the other hand, for $h \in \mathbb{C}[[x,y]]$, we define the intersection number
$$(f,h)_0=(C,C_h)_0:=dim_{\mathbb{C}} \frac{\mathbb{C} [[ x,y ]]}{(f,h)}=ord_t~h(x(t),y(t))$$
where $C_h$ is the Cartier divisor defined by $h$ and $\lbrace x(t)),y(t) \rbrace$ is as above.\\
The mapping $v_f:\frac {\mathbb{C} [[x,y]]}{(f)}\longrightarrow \mathbb{N}$, $h\longmapsto (f,h)_0$
defines a divisorial valuation. We define the semigroup of $C$ to be the semigroup of $v_f$ i.e
$\Gamma (C)=\Gamma (v_f)=\lbrace (f,h)_0 \in \mathbb{N},h \not \equiv 0 ~ mod(f) \rbrace$.\\
The following propositions and theorem from \cite{Z} characterize the structure of $\Gamma (C)$.

\begin{prop} There exists a unique sequence of $g+1$ positive integers
$(\bar{\beta_0}, \cdots ,\bar{\beta_g})$ such that:\\
$i) \bar{\beta_0}=\beta_0,$\\
$ii) \bar{\beta_i}=min\lbrace \Gamma (C)\backslash < \bar{\beta_0}, \cdots ,\overline{\beta}_{i-1}>\rbrace,1\leq i\leq g,$\\
$iii)\Gamma (C)=< \bar{\beta_0}, \cdots ,\bar{\beta_{g}}>,$\\
where for $i=1, \cdots ,g+1$,$ < \bar{\beta_0}, \cdots ,\overline{\beta}_{i-1}>$ is the semigroup generated by  $\bar{\beta_0}, \cdots ,\overline{\beta}_{i-1}$. By convention, we set $\bar{\beta}_{g+1}=+ \infty$.
\end{prop}
 
\begin{prop}The sequence $(\bar{\beta_0}, \cdots ,\bar{\beta_g})$ verifies:\\
$i)e_i=gcd(\bar{\beta_0}, \cdots ,\bar{\beta_{i}}),0 \leq i\leq g,$\\
$ii)\bar{\beta_0}=\beta_0$,$\bar{\beta_1}=\beta_1$ and $\bar{\beta_i}=\beta_i + \sum_{k=1}^{i-1}\frac{e_{k-1}-{e_k}}{e_{i-1}}\beta_k$,$i=2, \cdots ,g$.\\
$iii)n_i\bar{\beta_i}<\overline{\beta}_{i+1},1 \leq i\leq g-1$
\end{prop}
\begin{theo} The sequence $(\bar{\beta_0}, \cdots ,\bar{\beta_g})$ and the sequence $(\beta_0, \cdots .,\beta_g)$ are equivalent data
They determine and are determined by the topological type of $C$.
\end{theo}

Then from \cite{A} or \cite{Sp}, we can choose a system of approximate roots (or a minimal generating sequence)  $\{x_0, \cdots ,x_{g+1}\}$ of the divisorial valuation $v_f$.
We set $x=x_0, y=x_1$; for $i=2, \cdots ,g+1, x_i \in \mathbb{C}[[x,y]]$ is irreducible; for $1\leq i\leq g,$ the analytically irreducible
curve $C_i=\{x_i=0\}$ has $i-1$ Puiseux exponents and maximal contact with $C$ and $C_{g+1}=C$. This sequence also verifies\\
$i)$ $v_f(x_i)=\bar{\beta_i}$, $0\leq i \leq g,$\\
$ii) \Gamma (C_i)=<\frac{\bar{\beta_0}}{e_{i-1}}, \cdots ,\frac{\bar{\beta}_{i-1}}{e_{i-1}}>$ and the Puiseux sequence of
$C_i$ is $(\frac{\beta_1}{e_{i-1}}, \cdots ,\frac {\beta_{i-1}}{e_{i-1}})$,$2\leq i \leq g+1$.\\
$iii)$ for $1\leq i\leq g$, there exists a unique system of nonnegative integers $b_{ij}$,
$0\leq j<i$ such that for $1\leq j<i$, $b_{ij}<n_j$  and $n_i \bar{\beta_{i}}=\Sigma_{0\leq j<i}b_{ij}\bar{\beta_{j}}$.
And for $0\leq i\leq g$, one can choose $x_i$ such that they satisfy identities of the form

$$ x_{i+1}=x_i^{n_i}-c_ix_0^{b_{i0}}  \cdots  x_{i-1}^{b_{i(i-1)}}- \sum_{\gamma =(\gamma_0, \cdots ,\gamma_i)} c_{i,\gamma }x_0^{\gamma_0} \cdots x_i^{\gamma_i},(\star )$$
with ,$0\leq \gamma_j <n_j$, for $1\leq j<i$, and $\Sigma_{j}\gamma_j\bar{\beta_{j}}>n_i\bar{\beta_{i}}$
and with $c_{i,\gamma },c_i \in \mathbb{C}$ and $c_i \neq 0$.
These last equations $(\star )$ let us realize $C$ as a complete intersection in $\mathbb{C}^{g+1}=Spec ~\mathbb{C}~[[x_0, \cdots ,x_g]]$ defined by
the equations
$$ f_i=x_{i+1}-(x_i^{n_i}-c_ix_0^{b_{i0}} \cdots x_{i-1}^{b_{i(i-1)}}- \sum_{\gamma =(\gamma_0, \cdots ,\gamma_i)} c_{i,\gamma }x_0^{\gamma_0} \cdots x_i^{\gamma_i})$$
for $1\leq i\leq g$, with  $x_{g+1}=0$ by convention.\\

Let $h \in \mathbb{C}[[ x,y ]]$ be a $y$-regular irreducible power series with multiplicity
$p=ord_yh(0,y)$. Let $y(x^{\frac{1}{\beta_0}})$ and $z(x^{\frac{1}{p}})$ be  respectively roots
of $f$ and $g$ as in $(1)$. We call contact order of $f$ and $g$ in their Puiseux series the following rational
number $$o_f(h):=max\lbrace ord_x(y(wx^{\frac{1}{\beta_0}})-z(\lambda x^{\frac{1}{p}})); w^{\beta_0}=1,\lambda^{p}=1 \}=$$ $$max\lbrace ord_x(y(wx^{\frac{1}{\beta_0}})-z(x^{\frac{1}{p}}); w^{\beta_0}=1 \}=$$ 
$$max\lbrace ord_x(y(x^{\frac{1}{\beta_0}})-z(\lambda x^{\frac{1}{p}});\lambda^{p}=1 \}=o_h(f).$$
The following formula is from \cite{Me}, see also \cite{GP} .
\begin{prop}
Assume that $f$ and $h$ are as above; let $(\beta_1, \cdots ,\beta_g)$ the sequence of Puiseux exponents of $f$ and let $i \leq g+1$ be
the smallest strictly positive integer such that $o_f(h)\leq \frac{\beta_i}{\beta_0}$. Then
$$\frac{(f,h)_0}{p}=\sum_{k=1}^{i-1}\frac{e_{k-1}-{e_k}}{\beta_0}\beta_k +e_{i-1}o_f(h)$$
\end{prop}
\begin{cor}\cite{GP}
Let $i>0$ be an integer.Then  $o_f(h)\leq \frac{\beta_i}{\beta_0}$ iff $ \frac{(f,h)_0}{p}\leq e_{i-1}\frac{\bar{\beta_i}}{\beta_0}$.
Moreover  $o_f(h)= \frac{\beta_i}{\beta_0}$ iff $\frac{(f,h)_0}{p} =e_{i-1}\frac{\bar{\beta_i}}{\beta_0}$. In particular
$o_f(x_i)= \frac{\beta_i}{\beta_0},1\leq i\leq g.$

\end{cor}

\section{Jet schemes of complex branches}
We keep the notations of sections 2 and 3. We consider a curve $C \subset \mathbb{C}^2$ with a branch  of multiplicity $\beta_0>1$ at $0$, defined by $f$.
Note that in suitable coordinates we can write
$$f(x_0,x_1)=(x_1^{n_1}-cx_0^{m_1})^{e_1}+\sum_{a\beta_0+b\beta_1> \beta_0 \beta_1}c_{ab}x_0^ax_1^b;c \in \mathbb{C}^{\star}~ and ~c_{ab} \in \mathbb{C}. ~~~(\diamond )$$
We look for the irreducible components of $C_m^0:=(\pi_m^{-1}(0))$ for every $m \in \mathbb{N}$, where 
$\pi_m:C_m \rightarrow C $ is the canonical projection. Let $J_m^0$ be the radical of the ideal defining $(\pi_m^{-1}(0))$ in
$\mathbb{C}^2_m$.\\
In the sequel, we will denote the integral part of a rational number $r$ by $[r]$.
\begin{prop}\label{p1}
For $0<m<n_1\bar{\beta_1}$, we have that
$$(C_m^0)_{red}=(\pi_m^{-1}(0))_{red}=Spec~ \frac{\mathbb{C}[x_0^{(0)}, \cdots ,x_0^{(m)},x_1^{(0)}, \cdots ,x_1^{(m)}]}{(x_0^{(0)}, \cdots ,  x_0^{([\frac{m}{\beta_1}])},x_1^{(0)}, \cdots ,
x_1^{([\frac{m}{\beta_0}])})},$$
and
$$(C_{n_1\beta_1}^0)_{red}=(\pi_{n_1\beta_1} ^{-1}(0))_{red}=Spec \frac{\mathbb{C}[x_0^{(0)}, \cdots ,x_0^{(n_1\beta_1)},x_1^{(0)}, \cdots ,x_1^{(n_1\beta_1)}]}
{(x_0^{(0)}, \cdots ,x_0^{(n_1-1)},x_1^{(0)}, \cdots ,
x_1^{(m_1-1)},{x_1^{(m_1)}}^{n_1}-c{x_0^{(n_1)}}^{m_1})}.$$
\end{prop}
\begin{dem}
We write $f=\Sigma_{(a,b )}c_{ab}f_{ab}$ where $(a,b)\in \mathbb{N}^2,~f_{ab}=x_0^ax_1^b,~c_{ab} \in \mathbb{C}$
and $a\beta_0+b\beta_1\geq \beta_0 \beta_1$(the segment $[(0,\beta_0)(\beta_1,0)]$ is the Newton Polygon of $f$).
Let $supp(f)=\{(a,b)\in \mathbb{N}^2;c_{ab}\neq 0\}$.\\
For $0<m<n_1\beta_1$, the proof is by induction on $m$. For $m=1$,we have that $$F^{(1)}=\Sigma_{(a,b)\in supp(f)}c_{ab}F_{ab}^{(1)}$$
where $(F^{(0)}, \cdots ,F^{(i)})$ (resp.$(F_{ab}^{(0)}, \cdots ,F_{ab}^{(i)})$) is the ideal defining the $i$-th jet scheme $C_i$ of $C$(resp.
$C^{ab}_i$ the $i$-th jet scheme of $C^{ab}=\{f_{ab}=0\}$) in $\mathbb{C}^2_i$ .Then we have
$$ F_{ab}^{(1)}=\sum_{\sum i_k=1}x_0^{(i_1)} \cdots x_0^{(i_{a})}x_1^{(i_{a+1})} \cdots x_1^{(i_{a+b})}$$
where $\beta_1(a+b) \geq a\beta_0 +b\beta_1\geq \beta_0 \beta_1$ so $a+b \geq \beta_0>1$.
Then for every $(a,b)\in supp(f)$ and every $(i_1, \cdots ,i_a, \cdots ,i_{a+b})\in \mathbb{N}^{a+b}$ such that $\sum_{k=1}^{a+b}i_k=1$
there exists $1\leq k \leq a+b$ such that $i_{k}=0$, this means that $F_{ab}^{(1)}\in (x_0^{(0)},x_1^{(0)})$ and since we are looking over the origin, we
have that $(x_0^{(0)},x_1^{(0)})\subseteq J_1^0$ therefore
$(\pi_1^{-1}(0))_{red}=Spec \frac{\mathbb{C}[x_0^{(0)},x_0^{(1)},x_1^{(0)},x_1^{(1)}]}{(x_0^{(0)},x_1^{(0)})}$(In fact this is nothing but the Zariski tangent space of of $C$ at $0$).
\\ Suppose that the lemma holds until $m-1$ i.e.
$$(\pi_{m-1}^{-1}(0))_{red}=Spec \frac{\mathbb{C}[x_0^{(0)}, \cdots ,x_0^{(m-1)},x_1^{(0)}, \cdots ,x_1^{(m-1)}]}{(x_0^{(0)}, \cdots ,x_0^{([\frac{m-1}{\beta_1}])},x_1^{(0)}, \cdots ,
x_1^{([\frac{m-1}{\beta_0}])})}.$$
\underline{First case}:If $[\frac{m-1}{\beta_1}]=[\frac{m}{\beta_1}]$ and $[\frac{m-1}{\beta_0}]=[\frac{m}{\beta_0}]$.
We have $$F^{(m)}=\sum_{(a,b)\in supp(f)}c_{ab}\sum_{\sum i_k=m}x_0^{(i_1)} \cdots x_0^{(i_{a})}x_1^{(i_{a+1})} \cdots x_1^{(i_{a+b})}$$
Let $(a,b)\in supp(f)$; if for every $k=1, \cdots ,a,$ we had $i_k\geq [\frac{m}{\beta_1}]+1,$ and for every $k=a+1, \cdots ,a+b,$ we had $i_k\geq [\frac{m}{\beta_0}]+1,$
then $$m \geq a([\frac{m}{\beta_1}]+1)+b([\frac{m}{\beta_0}]+1)>\frac{m}{\beta_1}a +\frac{m}{\beta_0}b=m \frac{a\beta_0+b\beta_1}{\beta_0 \beta_1}\geq m.$$
The contradiction means that there exists $1\leq k\leq a$ such that $i_k\leq [\frac{m}{\beta_1}]$ or there exists 
$a+1\leq k\leq a+b$
such that $i_k \leq [\frac{m}{\beta_0}]$.
So $F^{(m)}$ lies in the ideal generated by  $J_{m-1}^0$ in $\mathbb{C}[x_0^{(0)}, \cdots ,x_0^{(m)},x_1^{(0)}, \cdots ,x_1^{(m)}]$ and $J_m^0=J^0_{m-1}.\mathbb{C}[x_0^{(0)}, \cdots ,x_0^{(m)},x_1^{(0)}, \cdots ,x_1^{(m)}].$
\\
\underline{Second case}:If $[\frac{m-1}{\beta_1}]=[\frac{m}{\beta_1}]$ and $[\frac{m-1}{\beta_0}]+1=[\frac{m}{\beta_0}]$
$(i.e.~\beta_0$ divides $m$). We have that
$$F^{(m)}=F_{0\beta_0}^{(m)}+ \sum_{(a,b)\in supp(f);(a,b)\neq (0,\beta_0) }F_{ab}^{(m)},~~(\star \star )$$
where
$$F_{0\beta_0}^{(m)}={\sum_{\sum i_k=m}}x_1^{(i_1)} \cdots x_1^{(i_{\beta_0})}={x_1^{(\frac{m}{\beta_0})}}^{\beta_0} +
\sum_{\sum i_k=m;(i_1, \cdots ,i_{\beta_0}) \neq (\frac{m}{\beta_0}, \cdots ,\frac{m}{\beta_0})}x_1^{(i_{1})} \cdots x_1^{(i_{\beta_0 })};$$
but $\sum i_k=m$ and $(i_1, \cdots ,i_{\beta_0}) \neq (\frac{m}{\beta_0}, \cdots ,\frac{m}{\beta_0})$ implies that there exists $1\leq k\leq \beta_0$
such that $i_k<\frac{m}{\beta_0},$ so
$$\sum_{\sum i_k=m;(i_1, \cdots ,i_{\beta_0}) \neq (\frac{m}{\beta_0}, \cdots ,\frac{m}{\beta_0})}x_1^{(i_{1})} \cdots x_1^{(i_{\beta_0 })} \in J_{m-1}^0.\mathbb{C}[x_0^{(0)}, \cdots ,x_0^{(m)},x_1^{(0)}, \cdots ,x_1^{(m)}].$$
For the same reason as above, we have that $$\sum_{(a,b)\in supp(f);(a,b)\neq (0,\beta_0) }F_{ab}^{(m)} \in J_{m-1}^0.\mathbb{C}[x_0^{(0)}, \cdots ,x_0^{(m)},x_1^{(0)}, \cdots ,x_1^{(m)}].$$
From $(\star \star )$ we deduce that ${x_1^{(\frac{m}{\beta_0})}} \in J_m^0$ and \\ $F^{(m)} \in (x_0^{(0)}, \cdots ,x_0^{([\frac{m}{\beta_1}])},x_1^{(0)}, \cdots ,
x_1^{(\frac{m}{\beta_0})}).$ Then $J_m^0=(x_0^{(0)}, \cdots ,x_0^{([\frac{m}{\beta_1}])},x_1^{(0)}, \cdots ,
x_1^{(\frac{m}{\beta_0})})$.\\
The third case i.e. if $[\frac{m-1}{\beta_1}]+1=[\frac{m}{\beta_1}]$ and $[\frac{m-1}{\beta_0}]=[\frac{m}{\beta_0}]$ is discussed as the second one.
Note that these are the only three possible cases since $m < n_1 \beta_1=lcm(\beta_0,\beta_1)$(here $lcm$ stands for the least common multiple).
\\
For $m=n_1\beta_1$, we have that $F^{(m)}$ is the coefficient of $t^m$ in the expansion of
$$f(x_0^{(0)}+x_0^{(1)}t+  \cdots  + x_0^{(m)}t^m,x_1^{(0)}+x_1^{(1)}t+  \cdots  + x_1^{(m)}t^m).$$
But since we are interested in the radical of the ideal defining the $m$-th
jet scheme, and we have found that $x_0^{(0)}, \cdots ,x_0^{(n_1-1)},x_1^{(0)}, \cdots ,
x_1^{(m_1-1)} \in J_{m-1}^0\subseteq J_{m}^0$, we can annihilate $x_0^{(0)}, \cdots ,x_0^{(n_1-1)},x_1^{(0)}, \cdots ,
x_1^{(m_1-1)}$ in the above expansion. Using   $(\diamond )$, we see that the coefficient
of $t^m$ is ${(x_1^{(m_1)}}^{n_1}-c{x_0^{(n_1)}}^{m_1})^{e_1}.$
\end{dem} 
$~~~$In the sequel if  $A$ is a ring ,  $I\subseteq A$ an ideal and $f \in A$, we denote by $V(I)$ the subvariety of $Spec~A$ defined by $I$
and by $D(f)$ the open set $\{f\neq 0\}$ in $Spec A ~i.e.~ D(f)=Spec A_f$.\\
The proof of the following corollary is analogous to that of proposition $4.1.$
\begin{cor}\label{pn}
Let $m \in \mathbb{N}$; let $k\geq 1$ be such that $m=kn_1\beta_1+i;1\leq i\leq n_1\beta_1$.
Then if $i<n_1\beta_1$, we have that $$Cont^{>kn_1}(x_0)_m=(\pi_{m,kn_1\beta_1} ^{-1}(V(x_0^{(0)}, \cdots ,x_0^{(kn_1)})))_{red}=$$
$$Spec ~\frac{k[x_0^{(0)}, \cdots ,x_0^{(m)},x_1^{(0)}, \cdots ,x_1^{(m)}]}
{(x_0^{(0)}, \cdots ,x_0^{(kn_1)}, \cdots ,x_0^{(kn_1+[\frac{i}{\beta_1}])},x_1^{(0)}, \cdots ,
x_1^{(km_1)}, \cdots ,x_1^{(km_1+[\frac{i}{\beta_0}])})}  $$
and if $i=n_1\beta_1$
$$ (\pi_{m,kn_1\beta_1} ^{-1}(V(x_0^{(0)}, \cdots ,x_0^{(kn_1)})))_{red}=$$
$$ Spec ~\frac{k[x_0^{(0)}, \cdots ,x_0^{(m)},x_1^{(0)}, \cdots ,x_1^{(m)}]}
{(x_0^{(0)}, \cdots ,x_0^{((k+1)n_1-1)},x_1^{(0)}, \cdots ,
x_1^{((k+1)m_1-1)},{x_1^{((k+1)m_1)}}^{n_1}-c{x_0^{((k+1)n_1)}}^{m_1})}.$$
\end{cor}

We now consider the case of a plane branch with one Puiseux exponent.
\begin{lem}
Let $C$ be a plane branch with one Puiseux exponent. Let $m,k \in \mathbb{N}$,
such that $k \neq 0 $ and $m\geq kn_1\beta_1+1$, and let $\pi_{m,kn_1\beta_1}:C_m \rightarrow C_{kn_1\beta_1}$ be the canonical projection. Then
$$C_m^k:=\pi_{m,kn_1\beta_1} ^{-1}(V(x_0^{(0)}, \cdots ,x_0^{(kn_1-1)})\cap D(x_0^{(kn_1)}))_{red}$$
is irreducible of codimension $k(m_1+n_1)+1 +(m-kn_1\beta_1)$ in $\mathbb{C}^2_{m}.$

\end{lem}
\begin{dem}
First note that since $e_1=1$, we have $m_1=\frac{\beta_1}{e_1}=\beta_1$.Let $I_m^{0k}$ be the ideal defining  
$C^k_m$ in $\mathbb{C}^2_{m} \cap D(x_0^{(kn_1)})$.Since $m \geq kn_1\beta_1$, by corollary $4.2$, $x_1^{(0)}, \cdots ,x_1^{(km_1-1)} \in
 I_m^{0k}$.So $I_m^{0k}$ is the radical of the ideal $I_m^{*0k}:=(x_0^{(0)}, \cdots ,x_0^{(kn_1-1)},x_1^{(0)}, \cdots ,x_1^{(km_1-1)},F^{(0)}, \cdots ,F^{(m)}).$ Now it follows from $\diamond$
and proposition $2.5$
that $$F^{(l)}\in (x_0^{(0)}, \cdots ,x_0^{(kn_1-1)},x_1^{(0)}, \cdots ,x_1^{(km_1-1)})~~for~~0\leq l< kn_1m_1,$$
$$F^{(kn_1m_1)}\equiv {x_1^{(km_1)}}^{n_1}-c{x_0^{(kn_1)}}^{m_1}~mod~(x_0^{(0)}, \cdots ,x_0^{(kn_1-1)},x_1^{(0)}, \cdots ,x_1^{(km_1-1)}),$$
$$F^{(kn_1m_1+l)}\equiv n_1{x_1^{(km_1)}}^{n_1-1}x_1^{(km_1+l)}-m_1c{x_0^{(kn_1)}}^{m_1-1}x_0^{(kn_1+l)}$$
$$+H_l(x_0^{(0)}, \cdots ,x_0^{(kn_1+l-1)},x_1^{(0)}, \cdots ,x_1^{(km_1+l-1)})~mod~(x_0^{(0)}, \cdots ,x_0^{(kn_1-1)},x_1^{(0)}, \cdots ,x_1^{(km_1-1)}),$$
for $1\leq l \leq m-kn_1m_1.$\\
This implies that 
$I_m^{*0k}:=(x_0^{(0)}, \cdots ,x_0^{(kn_1-1)},x_1^{(0)}, \cdots ,x_1^{(km_1-1)},F^{(kn_1m_1)}, \cdots ,F^{(m)}).$ Moreover the subscheme of
$\mathbb{C}^2_m \cap D(x_0^{(kn_1)})$ defined by  $I_m^{*0k}$ is isomorphic to the product of $\mathbb{C}^*$($\mathbb{C}^*$ is isomorphic to the regular locus of ${x_1^{(km_1)}}^{n_1}-c{x_0^{(kn_1)}}^{m_1}$)
by an affine space and its codimension is $k(m_1+n_1)+1 +(m-kn_1m_1)$; so it is reduced and irreducible, and it is nothing but
$C_m^k$, or equivalently $I_m^{0k}=I_m^{*0k}.$

\end{dem}

\begin{cor} Let $C$ be a plane branch with one Puiseux exponent. Let $m \in \mathbb{N}$,$m\neq 0$. let $q \in \mathbb{N}$ be
such that $m=qn_1\beta_1+i;0<i\leq n_1\beta_1$. Then $C_m^0=\pi_m^{-1}(0)$ has $q+1$ irreducible components which are:
$$C_{mkI}=\overline{C^k_m}  , 1\leq k\leq q,$$
$$and ~~B_m=Cont^{>qn_1}(x)_m=\pi_{m,qn_1\beta_1} ^{-1}(V(x_0^{(0)}, \cdots ,x_0^{(qn_1)})).$$
We have that $$codim(C_{mkI},\mathbb{C}^2_{m})=k(m_1+n_1)+1 +(m-kn_1m_1)$$ and
$$codim(B_m,\mathbb{C}^2_{m})=q(m_1+n_1)+ [\frac{i}{\beta_0}]+[\frac{i}{\beta_1}]+2=[\frac{m}{\beta_0}]+[\frac{m}{\beta_1}]+2 ~~if~~ i<n_1\beta_1$$
$$ codim(B_m,\mathbb{C}^2_{m})=(q+1)(m_1+n_1)+1 ~~if~~ i=n_1\beta_1.$$
\end{cor}
\begin{dem}
The codimensions and the irreducibility of $B_m$ and $C_{mkI}$  follow from  corollary 4.2 and lemma 4.3.  
  This shows that if $1\leq k<k'\leq q,
codim(C_{mk'I},\mathbb{C}^2_{m})<codim(C_{mkI},\mathbb{C}^2_{m})$ then $C_{mk'I} \nsubseteq C_{mkI}$.
On the other hand, since
$C_{mk'I} \subseteq  V(x_0^{(kn_1)})$ and $C_{mkI} \not\subseteq  V(x_0^{(kn_1)})$, we have that $C_{mkI} \nsubseteq C_{mk'I}$.
This also shows that $dim ~B_m \geq dim ~C_{mkI}$ for $1\leq k\leq q$, therefore
$B_m \not\subseteq C_{mkI}, 1\leq k\leq q$.But  $C_{mkI} \not\subseteq  B_m $ because $B_m \subseteq V(x_0^{(qn_1)})$ and
$C_{mkI} \not\subseteq V(x_0^{(qn_1)})$ for $1\leq k\leq q$. We thus have that $C_{mkI} \not\subseteq  B^m $ and $B^m \not\subseteq C_{mkI}$.
We conclude the corollary from the fact that by construction
$C_m^{0}=\cup_{k=1}^q C_{mkI} \cup B_m$.
\end{dem}
\\
~~~To understand the general case, i.e. to find the irreducible components of $C_m^{0}$ where $C$  has a branch with $g$ Puiseux exponents at $0$ , since for $kn_1\bar{\beta_{1}} < m \leq (k+1)n_1\bar{\beta_{1}}, m,k \in \mathbb{N} $ we know by corollary $4.2$ the structure of the $m$-jets that project to
$V(x_0^{(0)}, \cdots ,x_0^{(kn_1)})\cap C^{0}_{kn_1\beta_1}$, we search to understand for $m>kn_1\beta_1$ the $m$-jets that projects to \\
$V(x_0^{(0)}, \cdots ,x_0^{(kn_1-1)})\cap D(x_0^{(kn_1)})$, i.e. $C_m^k:=\pi_{m,kn_1\bar{\beta_{1}}}^{-1}(V(x_0^{(0)}, \cdots ,x_0^{(kn_1-1)})\cap D(x_0^{(kn_1)}))_{red}$. \\
Let $m, k \in  \mathbb{N}$ be such that $m\geq kn_1\beta_1$. Let $j=max\{l,n_2 \cdots n_{l-1}$ divides $k \}$(we set $j=2$ if
the greatest common divisor $(k,n_2)=1$ or if $g=1$).
Set $\kappa$ such that $k=\kappa n_2 \cdots n_{j-1}$, then we have $kn_1=
\kappa \frac{\beta_0}{n_j \cdots n_g}$.

\begin{prop}
Let $2 \leq j\leq g+1;$
for $i=2,..,g$,  
and $kn_1\bar{\beta_{1}}<m < \kappa e_{i-1}\frac{\bar{\beta_i}}{e_{j-1}}, $
we have that 
$$C_m^k=\bar{\pi}_{m, [\frac{m}{n_i \cdots n_g}]}^{-1}(C^k_{i,[\frac{m}{n_i \cdots n_g}]}), $$
where $\bar{\pi}_{m, [\frac{m}{n_i \cdots n_g}]}:\mathbb{C}^2_m\longrightarrow \mathbb{C}^2_{[\frac{m}{n_i \cdots n_g}]}$ is the canonical map. 
For $j<g+1$ and  $m \geq\kappa \bar{\beta_{j}}$,we have that $$C_m^k=\emptyset $$
\end{prop}

\begin{dem}
 Let $\phi  \in C_m^k$. Let $\tilde{\phi}:Spec\mathbb{C}[[t]]\longrightarrow (\mathbb{C}^2,0)$ be such that that lifts $\phi=\tilde{\phi}$ mod $ t^{m+1}.$ Let $\tilde{f} \in \mathbb{C}[[x,y]]$ be a function that defines the branch $\tilde{C}$ image of $\tilde{\phi}.$
we may assume that the map $Spec\mathbb{C}[[t]]\longrightarrow \tilde{C}$ induced by $\tilde{\phi}$ is the normalization of $\tilde{C}.$
Since $ord_t x_0 \circ \tilde{\phi}=kn_1, ord_tx_1\circ \tilde{\phi}=km_1, (ord_tx_0\circ \tilde{\phi}=kn_1)$ the multiplicity 
$m(\tilde{f})$ of $\tilde{C}$ at the origin is $ord_{x_1}\tilde{f} (0,x_1)=kn_1=\kappa \frac{\beta_{0}}{n_j \cdots n_g}.$\\
 \underline{Claim}: If $(f,\tilde{f})_{0} <\kappa e_{i-1}\frac{\bar{\beta_i}}{e_{j-1}} $ then  $(f,\tilde{f})_{0}=n_i \cdots n_g(x_i,\tilde{f})_0$.\\
Indeed, we have that $\frac{(f,\tilde{f})_{0}}{ord_y \tilde{f}(0,y)}<  e_{i-1}\frac{\bar{\beta_i}}{\beta_0}$,
therefore by corollary $3.5$ we have that $$o_f(\tilde{f}) < \frac{\beta_i}{\beta_0}=o_f(x_i).$$
Let $y(x^{\frac{1}{\beta_0}})$, $z(x^{\frac{1}{n_1 \cdots n_{i-1}}})$ and $u(x^{\frac{1}{m(\tilde{f})}})$
be respectively Puiseux-roots of $f$,$x_i$ and $\tilde{f}$. There exist $w$, $\lambda \in \mathbb{C}$ such that
$w^{\frac{\beta_0}{n_i \cdots n_g}}=1,\lambda^{m(\tilde{f})}=1$ and 
$$o_f(\tilde{f})=ord_x(u(\lambda x^{\frac{1}{m(\tilde{f})}})-y(x^{\frac{1}{\beta_0}}))$$ and
$$o_f(x_i)=ord_x(y(x^{\frac{1}{\beta_0}})-z(w x^{\frac{1}{n_1 \cdots n_{i-1}}})).$$
Since $o_f(\tilde{f})<o_f(x_i)$,  we have that 
$$o_f(\tilde{f})=ord_x(u(\lambda x^{\frac{1}{m(\tilde{f})}})-y(x^{\frac{1}{\beta_0}})+y(x^{\frac{1}{\beta_0}})-z(w x^{\frac{1}{n_1 \cdots n_{i-1}}}))$$
$$=ord_x(u(\lambda x^{\frac{1}{m(\tilde{f})}})-z(w x^{\frac{1}{n_1 \cdots n_{i-1}}}))\leq o_{x_i}(\tilde{f}).$$

On the other hand, there exist  $\lambda$ and $\delta  \in \mathbb{C}$, such that 
$\lambda^{m(\tilde{f})}=1,\delta^{\beta_0}=1$ and such that 
$$o_{x_i}(\tilde{f})=ord_x(u(\lambda x^{\frac{1}{m(\tilde{f})}})-z(x^{\frac{1}{ n_1 \cdots n_{i-1}}}))$$
and 
$$o_f(x_i)=ord_x(y(\delta x^{\frac{1}{\beta_0}})-z( x^{\frac{1}{n_1 \cdots n_{i-1}}})).$$
We have then that
$$o_{x_i}(\tilde{f})=ord_x(u(\lambda x^{\frac{1}{m(\tilde{f})}})-y(\delta x^{\frac{1}{\beta_0}}) +
y(\delta x^{\frac{1}{\beta_0}}) - z(w x^{\frac{1}{n_1 \cdots n_{i-1}}})).$$
Now $$ord_x(u(\lambda x^{\frac{1}{m(\tilde{f})}})-y(\delta x^{\frac{1}{\beta_0}})) \leq o_{f}(\tilde{f})<o_f(x_i)=
ord_x(y(\delta x^{\frac{1}{\beta_0}}) - z(w x^{\frac{1}{n_1 \cdots n_{i-1}}})).$$
So
$$o_{x_i}(\tilde{f})=ord_x(u(\lambda x^{\frac{1}{m(\tilde{f})}})-y(\delta x^{\frac{1}{\beta_0}}))\leq o_f(\tilde{f}).$$

We conclude that $ o_f(\tilde{f})=o_{x_i}(\tilde{f})$, and 
since the sequence of Puiseux exponents of $C_i$ is $(\frac{\beta_0}{n_i \cdots n_g}, \cdots ,\frac{\beta_{i-1} }{n_i \cdots n_g})$,
applying proposition $3.4$ to $C$ and $C_i$, we find that $(f,\tilde{f})_{0}=n_i \cdots n_g(x_i,\tilde{f})_0$ and claim follows.\\
On the other hand  by the corollary $3.5$ applied to $f$ and $\tilde{f}$,$(f,\tilde{f})_{0} \geq \kappa e_{i-1}\frac{\bar{\beta_i}}{e_{j-1}}$ if and only if
 $o_{f}(\tilde{f}) \geq \frac{\beta_i}{\beta_0}=o_{x_i}(f)=o_f(x_i) $ so $o_{f}(\tilde{f}) \geq \frac{\beta_i}{\beta_0}$ if and only if $o_{x_i}(\tilde{f})\geq \frac{\beta_i}{\beta_0}$, therefore $ (x_i,\tilde{f})_{0} \geq \kappa \frac{\bar{\beta_i}}{e_{j-1}}.$
This proves the first assertion. \\The second assertion is a direct consequence of 
lemma $5.1$ in \cite{GP}.
\end{dem}
To further analyse the $C_m^k$'s, we realize, as in section 3, $C$ as a complete intersection in $\mathbb{C}^{g+1}=Spec~\mathbb{C}[x_0, \cdots ,x_g]$ defined by
the ideal $(f_1, \cdots ,f_g)$ where
$$f_i= x_{i+1}-(x_i^{n_i}-c_ix_0^{b_{i0}} \cdots x_{i-1}^{b_{i(i-1)}}- \sum_{\gamma =(\gamma_0, \cdots ,\gamma_i)} c_{i,\gamma }x_0^{\gamma_0} \cdots x_i^{\gamma_i})$$
for $1\leq i\leq g$ and $x_{g+1}=0$. This will let us see the  $C_m^k$'s as fibrations over some reduced scheme that we understand well.\\
We keep the notations above and let $I_m^0$ be the radical of the ideal defining $C_m^{0}$ in $\mathbb{C}_m^{g+1}$  and let
$I_m^{0k}$ be the ideal defining $ C_m^k=(V(I_m^0,x_0^{(0)}, \cdots ,x_0^{(kn_1-1)})\cap D(x_0^{(kn_1)}))_{red}$ in $D(x_0^{(kn_1)})$. 

\begin{lem}   

Let $k \not=0$, $j$ and  $\kappa$ as above. For $1\leq i<j \leq g$ (resp.$1\leq i<j-1= g$)
and for $\kappa n_i \cdots n_{j-1}\bar{\beta_i}\leq m < \kappa n_{i+1} \cdots n_{j-1}\overline{\beta}_{i+1}$, we have
$$I_m^{0k}= (x_0^{(0)}, \cdots ,x_0^{( \frac{\kappa\bar{\beta_0}}{n_j \cdots n_g}-1)},$$
$$x_l^{(0)}, \cdots ,x_l^{( \frac{\kappa\bar{\beta_l}}{n_j \cdots n_g}-1)},F_l^{(\kappa \frac{n_l\bar{\beta_l}}{n_j \cdots n_g})}, \cdots ,F_l^{(m)},
1 \leq l\leq i, $$
$$x_{i+1}^{(0)}, \cdots ,x_{i+1}^{([\frac{m}{n_{i+1} \cdots n_g}])}, $$
$$F_l^{(0)}, \cdots ,F_l^{(m)},i+1\leq l \leq g-1).$$

Moreover for $1 \leq l\leq i$, 
$$F_l^{(\kappa \frac{n_l\bar{\beta_l}}{n_j \cdots n_g})} \equiv -({x_l^{(\kappa \frac{\bar{\beta_l}}{n_j \cdots n_g})}}^{n_l}-
c_l{x_0^{(\kappa \frac{\bar{\beta_0}}{n_j \cdots n_g})}}^{b_{l0}} \cdots .~~{x_{l-1}^{(\kappa \frac{\overline{\beta}_{l-1}}{n_j \cdots n_g})}}^{b_{l(l-1)}})$$
$$mod ~((x_l^{(0)}, \cdots ,x_l^{(\kappa \frac{\bar{\beta_l}}{n_j \cdots n_g}-1)})_{0 \leq l\leq i},x_{i+1}^{(0)}, \cdots ,x_{i+1}^{([\frac{m}{n_{i+1} \cdots n_g}])}),$$
for $1 \leq l< i$ and $\kappa \frac{n_l\bar{\beta_l}}{n_j \cdots n_g} <n<\kappa \frac{\overline{\beta}_{l+1}}{n_j \cdots n_g}$(resp. $l=i$ and $\kappa \frac{n_i\bar{\beta_i}}{n_j \cdots n_g} <n \leq [\frac{m}{n_{i+1} \cdots n_g}]$)
$$F_l^{(n)}\equiv -(n_lx_l^{(\kappa \frac{\bar{\beta_l}}{n_j \cdots n_g})^{n_l-1}}
x_l^{(\kappa \frac{\bar{\beta_l}}{n_j \cdots n_g}+n-
\kappa \frac{n_l\bar{\beta_l}}{n_j \cdots n_g})}-$$
$$c_l\sum_{0\leq h \leq l-1}b_{lh}x_0^{(\kappa\frac{\bar{\beta_0}}{n_j \cdots n_g})^{b_{l0}}} \cdots x_h^{(\kappa\frac{\bar{\beta_h}}{n_j \cdots n_g})^{b_{lh}-1}}
x_h^{(\kappa\frac{\bar{\beta_h}}{n_j \cdots n_g}+n-
\kappa\frac{n_l\bar{\beta_l}}{n_j \cdots n_g})} \cdots x_{l-1}^{(\kappa
\frac{\overline{\beta_{l-1}}}{n_j \cdots n_g})^{b_{l(l-1)}}}+$$
$$H_l( \cdots ,x_h^{(\kappa\frac{\bar{\beta_h}}{n_j \cdots n_g}+n-\kappa \frac{n_l\bar{\beta_l}}{n_j \cdots n_g}-1)}, \cdots ) )$$
$$mod ~((x_l^{(0)}, \cdots ,x_l^{(\kappa \frac{\bar{\beta_l}}{n_j \cdots n_g}-1)})_{0 \leq l\leq i},x_{i+1}^{(0)}, \cdots ,x_{i+1}^{([\frac{m}{n_{i+1} \cdots n_g}])}),$$

for $1 \leq l< i$ and $\kappa \frac{\overline{\beta_{l+1}}}{n_j \cdots n_g} \leq n \leq m $(resp. $l=i$ and $[\frac{m}{n_{i+1} \cdots n_g}] <n \leq m$), or $i+1\leq l \leq g-1$
and $ 0 \leq n \leq m$,
$$ F_l^{(n)}=x_{l+1}^{(n)}+H_l(x_0^{(0)}, \cdots ,x_0^{(n)}, \cdots ,x_l^{(0)}, \cdots ,x_l^{(n)}).$$
For $i=j-1=g$ and $ m \geq \kappa n_g\bar{\beta_g}$,
$$I_m^{0k}= (x_0^{(0)}, \cdots ,x_0^{( \kappa\bar{\beta_0}-1)},$$
$$x_l^{(0)}, \cdots ,x_l^{(\kappa\bar{\beta_l}-1)},F_l^{(\kappa n_l\bar{\beta_l})}, \cdots ,F_l^{(m)}),1 \leq l \leq g,$$
where for $1 \leq l < g $ and $\kappa n_l \bar{\beta_l} \leq n \leq m$, the above formula for $F_l^{(n)}$ remains valid,
$$ F_g^{(\kappa n_g\bar{\beta_g})} \equiv 
-({x_g^{(\kappa \bar{\beta_g})^{n_g}}-
c_gx_0^{(\kappa \bar{\beta_0})^{b_{g0}}} \cdots .~~x_{g-1}^{(\kappa \overline{\beta_{g-1}})^{b_{g(g-1)}}}})$$
$$mod ~((x_l^{(0)}, \cdots ,x_l^{(\kappa\bar{\beta_l}-1)}))_{0 \leq l\leq g}$$
and for $\kappa n_g\bar{\beta_g}< n \leq m,$
$$F_g^{(n)}\equiv -(n_gx_g^{(\kappa \bar{\beta_g})^{n_g-1}}
x_g^{(\kappa \bar{\beta_g}+n-\kappa n_g\bar{\beta_g})}-$$
$$c_g\sum_{0\leq h \leq g-1}b_{g0}x_0^{(\kappa\bar{\beta_0})^{b_{gh}}} \cdots x_h^{(\kappa \bar{\beta_h})^{b_{gh}-1}}
x_h^{(\kappa\bar{\beta_h}+n-\kappa n_h\bar{\beta_h})} \cdots x_{g-1}^{(\kappa
\overline{\beta_{g-1}})^{b_{g(g-1)}}}+$$
$$H_g( \cdots ,x_h^{(\kappa\bar{\beta_h}+n-\kappa n_h\bar{\beta_h})}, \cdots ) )$$
$$mod ~((x_l^{(0)}, \cdots ,x_l^{(\kappa\bar{\beta_l}-1)}))_{0 \leq l\leq g}$$

\end{lem}

\begin{dem}
First assume that $\kappa n_i \cdots n_{j-1}\bar{\beta_i}\leq m < \kappa n_{i+1} \cdots n_{j-1}\bar{\beta}_{i+1}$ for $1\leq i<j\leq g$
(resp. $1 \leq i < j-1=g)$. By proposition $4.5$, we have that $C_m^k=\bar{\pi}_{m,[\frac{m}{n_{i+1} \cdots n_g}]}^{-1}(C_{i+1,[\frac{m}{n_{i+1} \cdots n_g}]}^k)$ where
 $\bar{\pi}_{m,[\frac{m}{n_{i+1} \cdots n_g}]}:\mathbb{C}^2_m \longrightarrow \mathbb{C}^2_{[\frac{m}{n_{i+1} \cdots n_g}]}$
is the canonical map. Now $\mathbb{C}^2=Spec ~\mathbb{C}[x_0,x_1] (resp.$ $C_{i+1}=V(x_{i+1}))$ is realized as the complete intersection in 
$\mathbb{C}^{g+1}=Spec ~\mathbb{C}[x_0, \cdots ,x_g]$ defined by the ideal $(f_1, \cdots ,f_{g-1})$(resp. $(f_1, \cdots ,f_{g-1},x_{i+1})$). So
since 
$m \geq kn_1\bar{\beta_1} , 
I_m^{0k}$ is the radical of the 
ideal $I_m^{*0k}=$ $$(x_0^{(0)}, \cdots ,x_0^{(kn_1-1)},x_1^{(0)}, \cdots ,
x_1^{(km_1-1)},F_1^{(0)}, \cdots ,F_1^{(m)},$$ $$ \cdots ,F_{g-1}^{(0)}, \cdots ,F_{g-1}^{(m)},x^{(0)}_{i+1}, \cdots ,x^{([\frac{m}{n_{i+1} \cdots n_g}])}_{i+1}).$$
We first observe that $F_1^{(n)}\equiv x_2^{(n)}~ mod~ (x_0^{(0)}, \cdots ,x_0^{(kn_1-1)},x_1^{(0)}, \cdots ,
x_1^{(km_1-1)})$ for $0\leq n <kn_1\bar{\beta_1}.$~Now since $\frac{m}{n_{2} \cdots n_g} \geq [\frac{m}{n_{2} \cdots n_g}] \geq kn_1m_1,$ we have
$$F_1^{(kn_1m_1)} \equiv-(x_1^{(km_1)^{n_1}}-c_1x_0^{(kn_1)^{m_1}})~$$
$$mod~(x_0^{(0)}, \cdots ,x_0^{(kn_1-1)},x_1^{(0)}, \cdots ,
x_1^{(km_1-1)},x_2^{(0)}, \cdots , x_2^{([\frac{m}{n_{2} \cdots n_g}])})$$ and 
$$F_1^{(n)} \equiv-(n_1x_1^{(km_1)^{n_1-1}}x_1^{(km_1+n-kn_1m_1)}-m_1c_1x_0^{(kn_1)^{m_1-1}}x_0^{(kn_1+n-kn_1m_1)})$$
$$+H_1(x_0^{(0)}, \cdots ,x_0^{(kn_1+n-kn_1m_1-1)},x_1^{(0)}, \cdots ,x_1^{(km_1+n-kn_1m_1-1)})$$
$$~mod~(x_0^{(0)}, \cdots ,x_0^{(kn_1-1)},x_1^{(0)}, \cdots ,
x_1^{(km_1-1)},x_2^{(0)}, \cdots , x_2^{([\frac{m}{n_{2} \cdots n_g}])})$$ for $kn_1\bar{\beta_1}<n\leq [\frac{m}{n_{2} \cdots n_g}]$.
Finally, for $l=1$ and $[\frac{m}{n_{2} \cdots n_g}] <n \leq m,$ or $2\leq l \leq g-1$ and $0 \leq n \leq m$, 
we have $$F_l^{(n)}=x^{(n)}_{l+1}+H_l(x^{(0)}_0, \cdots ,x_0^{(n)}, \cdots ,x_l^{(0)}, \cdots ,x_l^{(n)}).$$
As a consequence for  $i=1$, the subscheme of $\mathbb{C}^{g+1} \cap D(x_0^{(kn_1)})$ defined by $I_m^{*0k}$ is isomorphic to the product of $\mathbb{C}^*$ by an affine space , so it is reduced and irreducible and  $I_m^{*0k}= I_m^{0k}$ is a prime ideal 
in $\mathbb{C}[x_0^{(0)}, \cdots ,x_0^{(m)}, \cdots ,x_g^{(0)}, \cdots ,x_g^{(m)}]_{x_0^{(kn_1)}}$, generated by a regular sequence, i.e the proposition holds for $i=1$.\\
Assume that it holds for $i<j-1<g$(resp. $i<j-2=g-1).$
For $\kappa n_{i+1} \cdots n_{j-1}\overline{\beta}_{i+1}\leq m < \kappa n_{i+2} \cdots n_{j-1}\overline{\beta}_{i+2}$, the ideal  in $\mathbb{C}[x_0^{(0)}, \cdots ,x_0^{(m)}, \cdots ,x_g^{(0)}, \cdots ,x_g^{(m)}]_{x_0^{(kn_1)}}$ generated
by $I_{\kappa n_{i+1} \cdots n_{j-1}\overline{\beta_{i+1}}-1}^{0k}$ is contained in $I_m^{0k}.$ By the inductive hypothesis,
$x_l^{(0)}, \cdots ,x_l^{(\frac{\kappa\bar{\beta_l}}{n_j \cdots n_g}-1)} \in 
I^{0k}_{\kappa n_{i+1} \cdots n_{j-1}\overline{\beta}_{i+1}-1}$
for $l=1, \cdots ,i+1.$ So $I_m^{0k}$ is the radical of
$$I_m^{*0k}=(x_0^{(0)}, \cdots ,x_0^{(\frac{\kappa\bar{\beta_0}}{n_j \cdots n_g}-1)},$$
$$x_l^{(0)}, \cdots ,x_l^{(\frac{\kappa\bar{\beta_l}}{n_j \cdots n_g}-1)},F_l^{(0)}, \cdots ,F_l^{(m)},1\leq l \leq i+1,$$
$$  x_{i+2}^{(0)}, \cdots ,x_{i+2}^{([\frac{m}{n_{i+2} \cdots n_g}])},$$
$$ F_l^{(0)}, \cdots ,F_l^{(m)},i+2\leq l \leq g-1).$$
Now  for $0 \leq n < \frac{\kappa n_l\bar{\beta_l}}{n_j \cdots n_g}$,we have $$F_l^{(n)} \equiv x_{l+1}^{(n)} ~mod~(x_0^{(0)}, \cdots ,x_l^{(\frac{\kappa\bar{\beta_0}}{n_j \cdots n_g}-1)},x_l^{(0)}, \cdots ,x_l^{(\frac{\kappa\bar{\beta_l}}{n_j \cdots n_g}-1)},$$
$$1\leq l \leq i+1).$$
Here since $\overline{\beta}_{l+1}>n_l\bar{\beta_l},$ for $1 \leq l \leq i$ and $\frac{m}{n_{i+2} \cdots n_g} \geq [\frac{m}{n_{i+2} \cdots n_g}]
\geq \frac{\kappa n_{i+1}\overline{\beta}_{i+1}}{n_j \cdots n_g}$, we can delete $F_l^{(n)}$, $1 \leq l \leq i+1, 0 \leq n < \frac{\kappa n_l\bar{\beta_l}}{n_j \cdots n_g}$ from the above generators of $I_m^{*0k}$ without changing the generated ideal. The identities relative to the $F_l^{(n)}$ for $1 \leq l \leq i+1,  \frac{\kappa n_l\bar{\beta_l}}{n_j \cdots n_g} \leq n \leq m$ or
$i+2 \leq l \leq g-1$ and $0 \leq n \leq m$ follow immediately from  $(\diamond).$ So here again the 
subscheme of $\mathbb{C}^{g+1} \cap D(x_0^{(kn_1)})$ defined by $I_m^{*0k}$ is isomorphic to the product of $\mathbb{C}^*$ by an affine space , so it is reduced and irreducible and  $I_m^{*0k}= I_m^{0k}$ is a prime ideal 
in $\mathbb{C}[x_0^{(0)}, \cdots ,x_0^{(m)}, \cdots ,x_g^{(0)}, \cdots ,x_g^{(m)}]_{x_0^{(kn_1)}}$, generated by a regular sequence, i.e the proposition holds for $i+1$.\\
The case $i=j-1=g$ and $m \geq \kappa n_g\overline{\beta_g}$ follows by similar arguments.
\end{dem}
As an immediate consequence we get
\begin{prop} \label{dim}
Let $C$ be a plane branch with $g$ Puiseux exponents. Let $k \not=0, j$ and $\kappa$ as above.
 For $m\geq kn_1\beta_1$,  let $\pi_{m,kn_1\beta_1}:C_m \rightarrow C_{kn_1\beta_1}$ be the canonical projection and let
$C_m^k:=\pi_{m,kn_1\beta_1} ^{-1}(D(x_0^{(kn_1)})\cap V(x_0^{(0)}, \cdots ,x_0^{(kn_1-1)}))_{red}.$ Then for $1 \leq i <j \leq g$
(resp.$1 \leq i <j-1=g) $ and $\kappa n_i \cdots n_{j-1} \bar{\beta_i} \leq m < \kappa n_{i+1} \cdots n_{j-1} \overline{\beta}_{i+1}$,
$C_m^k$ is irreducible of codimension 
$$\frac{\kappa}{n_j \cdots n_g}( \bar{\beta_0} +\bar{\beta_1} + \sum_{l=1}^{i-1}(\overline{\beta}_{l+1}-n_l\overline{\beta_{l}}))
+([\frac{m}{n_{i+1} \cdots n_g}] -\frac{\kappa n_i\bar{\beta_i}}{n_j \cdots n_g})+1$$
in $\mathbb{C}^2_{m}$.\\
For $j\leq g$ and $m \geq \kappa\bar{\beta_j}$ (resp.$ j=g+1$ and $m \geq \kappa n_g\bar{\beta_g}$),
$$C_m^k=\emptyset$$
(resp. $C_m^k$ is of codimension  $$\kappa( \bar{\beta_0} +\bar{\beta_1} + \sum_{l=1}^{g-1}(\overline{\beta}_{l+1}-n_l\overline{\beta_{l}}))
+m -\kappa n_g\bar{\beta_g}+1) $$
in $\mathbb{C}^2_{m}$.

\end{prop}
For $k' \geq k$ and $m \geq k'n_1\beta_1$, we now compare codim($C_m^{k}$,$\mathbb{C}^2_{m}$) and codim($C_m^{k'}$,$\mathbb{C}^2_{m}$).
\begin{cor}\label{cd}
 For $k' \geq k \geq 1$ and $m\geq k'n_1\beta_1$, if $C_m^{k}$ and $C_m^{k'}$ are nonempty, we have
$$codim(C_m^{k'},\mathbb{C}^2_{m}) \leq codim(C_m^{k},\mathbb{C}^2_{m}).$$
\end{cor}
\begin{dem}
 Let $\gamma^k:[kn_1\beta_1, \infty [\longrightarrow [k(n_1+m_1), \infty [$ be the function given by 
$$\gamma^k(m)=\frac{k}{e_1}( \bar{\beta_0} +\bar{\beta_1} + \sum_{l=1}^{i-1}(\overline{\beta}_{l+1}-n_l\overline{\beta_{l}}))
+(\frac{m}{e_i} -\frac{k n_i\bar{\beta_i}}{e_1})+1$$
for $1\leq i <g $ and  $\frac{k\bar{\beta_i}}{n_2 \cdots n_{i-1}} \leq m < \frac{k\overline{\beta}_{i+1}}{n_2 \cdots n_{i}}$
and 
$$\gamma^k(m)=\frac{k}{e_1}( \bar{\beta_0} +\bar{\beta_1} + \sum_{l=1}^{g-1}(\overline{\beta}_{l+1}-n_l\overline{\beta}_l))
+(m -\frac{k n_g\bar{\beta_g}}{e_1})+1$$
for $i=g$ and $m \geq \frac{k\overline{\beta_{g}}}{n_2 \cdots n_{g-1}}.$
\\
In view of proposition \ref{dim} , we have that codim($C_m^{k}$,$\mathbb{C}^2_{m})=[\gamma^k(m)]$ for $k \equiv 0$ mod $n_2 \cdots n_{j-1}$
and $k \not \equiv 0$ mod $ n_2 \cdots n_{j}$ with $2\leq j \leq g$ and any integer $m \in [kn_1\beta_1,\frac{k\overline{\beta}_j}{n_2 \cdots n_{j-1}}[$
or for $k \equiv 0$ mod $n_2 \cdots n_g$and any integer $m\geq kn_1\beta_1$. Similarly we define $\gamma^{k'}:[k'n_1\beta_1, \infty [\longrightarrow [k'(n_1+m_1), \infty [$ by changing $k$ to $k'.$\\
Let $\Gamma^k(resp.\Gamma^{k'})$ be the graph of $\gamma^k$(resp $\gamma^{k'}$) in $\mathbb{R}^2$.Now let $\tau:\mathbb{R}^2\longrightarrow \mathbb{R}^2$ be defined by $\tau(a,b)=(a,b-1)$ and let $\lambda^{k'/k}:\mathbb{R}^2\longrightarrow \mathbb{R}^2$ be defined by $\lambda^{k'/k}(a,b)=\frac{k'}{k}(a,b).$ We note that $\tau(\Gamma^{k'})=\lambda^{k'/k}(\tau(\Gamma^{k}));$
we also note that the endpoints of $\tau(\Gamma^{k})$ and $\tau(\Gamma^{k'})$ lie on the line through $0$ with slope $\frac{\beta_0+\beta_1}{e_1n_1\beta_1}=
\frac{1}{e_1}\frac{n_1+m_1}{n_1m_1}<\frac{1}{e_1}$. Since $\frac{k'}{k} \geq 1$, the image of $\tau(\Gamma^{k})$ by $\lambda^{k'/k}$
lie on the subset of $\mathbb{R}^2$ whith boundary the union of $\tau(\Gamma^{k})$, of the segment joining its endpoint
$(kn_1\beta_1,\frac{k}{e_1}(\beta_0+\beta_1))$ to $(kn_1\beta_1,0)$ and of $[kn_1\beta_1,\infty[ ~\times~0.$ This implies that 
$\gamma^{k'}(m) \leq \gamma^{k}(m)$ for $m \geq k'n_1 \beta_1$ ,hence $ [\gamma^{k'}(m)] \leq [\gamma^{k}(m)]$ and the claim.
 
\end{dem}

\begin{theo}
Let $C$ be a plane branch with $g \geq 2$ Puiseux exponents. Let $m \in \mathbb{N}$. For $1\leq m < n_1\beta_1+e_1$,$C_m^0=Cont^{>0}(x_0)_m$ is irreducible.
For $qn_1\beta_1+e_1\leq m <(q+1)n_1\beta_1+e_1$,with $q\geq 1$ in $\mathbb{N}$, 
 the irreducible components of $C_m^0$ are :
$$C_{m\kappa I}=\overline{Cont^{\kappa\bar{\beta_0}}(x_0)_m}$$ 
for $1\leq \kappa$ and $\kappa\bar{\beta_0}\bar{\beta_1}+e_1 \leq m,$\\
 $$C^j_{m\kappa v}=\overline{Cont^{\frac{\kappa\bar{\beta_0}}{n_j \cdots n_g}}(x_0)_m} $$ 
for $j=2, \cdots ,g,$  $1\leq \kappa$ and $\kappa \not \equiv 0~ mod ~n_j$ and such that
$\kappa n_1 \cdots n_{j-1}\bar{\beta_1}+e_1\leq m< \kappa\bar{\beta_j}$,\\
 
$$B_m= Cont^{>n_1q}(x_0)_m.$$
\end{theo}

\begin{dem}
 We first observe that for any integer $k \not = 0$ and any $m \geq kn_1\beta_1$,
$$(C_m^0)_{red}=\cup_{1\leq h \leq k} C_m^h\cup Cont^{>kn_1}(x_0)_m $$
where $C_m^h:=Cont^{hn_1}(x_0)_m$ as above. Indeed , for $k=1$, we have that 
$(C_m^0)_{red} \subset V(x_0^{(0)}, \cdots ,x_0^{(n_1-1)})$ by proposition \ref{p1}.
Arguing by induction on $k$, we may assume that the claim holds for $m\geq(k-1)n_1\beta_1$.Now by corollary
\ref{pn}, we know that for $m\geq kn_1\beta_1$, $Cont^{>(k-1)n_1}(x_0)_m\subset  V(x_0^{(0)}, \cdots ,x_0^{(kn_1-1)}),$
hence the claim for $m\geq  kn_1\beta_1$.\\
We thus get that for $qn_1\beta_1+e_1\leq m <(q+1)n_1\beta_1+e_1,$
$$(C_m^0)_{red}=\cup_{1\leq k \leq q} C_m^k\cup Cont^{>qn_1}(x_0)_m .$$
By proposition \ref{dim},for $1 \leq k \leq q,C_m^k$ is either irreducible or empty. We first note that if $C_m^k\not=\emptyset,$ then 
$\overline{C^k_m}\not\subset Cont^{>qn_1}(x_0)_m.$ Similarly, if $1 \leq k<k'\leq q$ and if $C^k_m$ and $C^{k'}_m$ are nonempty,
then $\overline{C^{k}_m} \not\subset \overline{C^{k'}_m}.$ On the other hand by corollary  \ref{cd}, we have that
$codim(C_m^{k'},\mathbb{C}^2_{m}) \leq codim(C_m^{k},\mathbb{C}^2_{m})$. So $\overline{C^{k'}_m} \not\subset \overline{C^k_m}.$
Finally we will show that $Cont^{>qn_1}(x_0)_m$ $ \not\subset \overline{C^k_m}$ if $C^k_m \not= \emptyset$ for $1 \leq k \leq q$.
To do so, it is enough to check that $codim(C_m^{k},\mathbb{C}^2_{m}) \geq codim(Cont^{>qn_1}(x_0)_m,\mathbb{C}^2_{m}).$
For $m \in [qn_1\beta_1+e_1 ,(q+1)n_1\beta_1[$, we have 
$$\delta^q(m):=codim(Cont^{>qn_1}(x_0)_m,\mathbb{C}^2_{m})=2+q(n_1+m_1)+[\frac{m-qn_1\beta_1}{\beta_0}]+[\frac{m-qn_1\beta_1}{\beta_1}]$$
by corollary $4.2.$Let $\lambda^q:[qn_1\beta_1+e_1[ \longrightarrow [q(n_1+m_1),\infty[$ be the function given by  
$\lambda^q(m)=q(n_1+m_1)+\frac{m-qn_1\beta_1}{e_1}+1.$
For simplicity, set $i=m-qn_1\beta_1$.For any integer $i$ such that 
$e_1 \leq i <n_1\beta_1=n_1m_1e_1$, we have $1+[\frac{i}{n_1e_1}]+[\frac{i}{m_1e_1}]\leq [\frac{i}{e_1}].$
Indeed this is true for $i=e_1$ and it follows by induction on $i$ from the fact that for any pair of integers $(b,a),$
we have $[\frac{b+1}{a}]=[\frac{b}{a}]$ if and only if $b+1 \not \equiv 0 $ mod $a$ and $[\frac{b+1}{a}]=[\frac{b}{a}]+1$ otherwise,
since $i <n_1m_1e_1.$ So $\delta^q(m) \leq [\lambda^q(m)].$ \\
But in the proof of corollary \ref{cd}, we have checked that if $C_m^{k} \not=\emptyset,$ we have 
codim($C_m^{k}$,$\mathbb{C}^2_{m})=[\gamma^k(m)]$. We have also checked that for $q\geq k$ and $m\geq qn_1\beta_,$ 
$\gamma^k(m) \geq \gamma^q(m).$ Finally in view of the definitions of $\gamma^q$ and $\lambda^q,$ we have $\gamma^q(m) \geq \lambda^q(m),$ 
so $[\gamma^q(m)] \geq [\lambda^q(m)]\geq \delta^q(m).$ \\
For $m=(q+1)n_1\beta_1,$ we have $\delta^q(m)=(q+1)(n_1+m_1)+1$ by corollary \ref{pn}. For $m \in [(q+1)n_1\beta_1,(q+1)n_1\beta_1+e_1[,$
we have $Cont^{>qn_1}(x_0)_m=C_m^{q+1}\cup Cont^{>(q+1)n_1}(x_0)_m$ and $Cont^{>(q+1)n_1}(x_0)_m=V(x_0^{(0)}, \cdots ,x_0^{((q+1)n_1)},
x_1^{(0)}, \cdots ,x_1^{((q+1)m_1)})$ again by corollary \ref{pn}. If in addition we have $m<(q+1)\bar{\beta_2},$
then by proposition $4.5$ $C_m^{q+1}=V(x_0^{(0)}, \cdots ,x_0^{((q+1)n_1-1)},$ $x_1^{(0)}, \cdots ,x_1^{((q+1)m_1-1)},x_1^{((q+1)m_1)^{n_1}}-c_1x_0^{((q+1)n_1)^{m_1}})\cap 
D(x_0^{((q+1)n_1)},$ thus we have \\ $Cont^{>qn_1}(x_0)_m$ $=\overline{C_m^{q+1}}$ and $\delta^q(m)=(q+1)(n_1+m_1)+1.$
We have $(q+1)n_1\beta_1+e_1 \leq (q+1)
\overline{\beta}_2$ if $q+1  \geq  n_2,$ because $\overline{\beta}_2 -n_1\overline{\beta}_1 \equiv 0$ mod $(e_2)$ . If not , we may have $(q+1)\overline{\beta}_2<(q+1)n_1\beta_1+e_1$, so for 
$(q+1)\overline{\beta_2}\leq m<(q+1)n_1\beta_1+e_1,$ we have $C_m^{q+1}=\emptyset, Cont^{>qn_1}(x_0)_m=Cont^{>(q+1)n_1}(x_0)_m$
and $\delta^q(m)=(q+1)(n_1+m_1)+2.$\\
In both cases, for $m \in [(q+1)n_1\beta_1,(q+1)n_1\beta_1+e_1[,$ we have $\delta^q(m) \leq (q+1)(n_1+m_1)+2.$ Since
$[\lambda^q(m)]=q(n_1+m_1)+n_1m_1+1$, we conclude that $[\lambda^q(m)] \geq \delta^q(m)$,
so for $1 \leq k \leq q,$ if $C_m^k \not= \emptyset,$ we have $[\gamma^k(m)] \geq \delta^q(m).$ This proves that the irreducible components of $C_m^{0}$ are the $\overline{C_m^k}$ for $1 \leq k \leq q$ and $C_m^k \not= \emptyset$, and $Cont^{>qn_1}(x_0)_m,$
hence the claim in viewof the characterization of the nonempty $C_m^{k's}$'s given in proposition $4.5.$

\end{dem}

\begin{cor}
Under the assumption of theorem $4.9$, let $q_0+1=min\{\alpha \in \mathbb{N};\alpha(\overline{\beta}_2-n_1\overline{\beta}_1)\geq e_1\}$.
Then $0\leq q_0 <n_2.$
For $1\leq m<(q_0+1)n_1\beta_1+e_1, C_m^0$ is irreducible and we have $codim(C_m^0,\mathbb{C}^2_m)=$
 
$$2+[\frac{m}{\beta_0}]+[\frac{m}{\beta_1}] ~~for ~~0\leq q \leq q_0~~and ~~ qn_1\beta_1+e_1 \leq m <(q+1)n_1\beta_1~~ $$ $$ or ~~0\leq q \leq q_0~~ and ~~(q+1)\overline{\beta}_2\leq m< (q+1)n_1\beta_1+e_1.$$
$$1+[\frac{m}{\beta_0}]+[\frac{m}{\beta_1}] ~~for ~~0\leq q < q_0~~and ~~ (q+1)n_1\beta_1 \leq m <(q+1)\overline{\beta}_2~~ $$ $$ or ~~(q_0+1)n_1\beta_1\leq m< (q_0+1)n_1\beta_1+e_1.$$

For $q\geq q_0+1$ in $\mathbb{N}$ and  $qn_1\beta_1+e_1 \leq m <(q+1)n_1\beta_1+e_1,$ the number of irreducible components of 
  $C_m^0$  is:
 $$N(m)=q+1-\sum^g_{j=2}([\frac{m}{\bar{\beta_j}}]-[\frac{m}{n_j\bar{\beta_j}}])$$
and 
$codim(C_m^0,\mathbb{C}^2_m)=$
$$2+[\frac{m}{\beta_0}]+[\frac{m}{\beta_1}] ~~for~~   qn_1\beta_1+e_1 \leq m <(q+1)n_1\beta_1.~~$$
$$1+[\frac{m}{\beta_0}]+[\frac{m}{\beta_1}] ~~for~~ (q+1)n_1\beta_1 \leq m  <(q+1)n_1\beta_1+e_1. $$
\end{cor}
\begin{dem}
 We have already observed that $n_2(\overline{\beta}_2-n_1\overline{\beta}_1)\geq e_1$ because 
$\overline{\beta}_2-n_1\overline{\beta}_1 \equiv 0$ mod $(e_2),$ so $1\leq q_0+1\leq n_2$. \\
For $qn_1\beta_1+e_1 \leq m <(q+1)n_1\beta_1+e_1,$  with $q\geq 1$, we have seen in the proof of theorem $4.9$ that the irreducible components of $C_m^0$ are the $\overline{C_m^k}$ for $1\leq k \leq q$ and $C_m^k\not= \emptyset$ and $Cont^{qn_1}(x_0)_m$. We thus have to enumerate the empty $C_m^k$ for  $1\leq k \leq q.$ By proposition $4.5$, $C_m^k= \emptyset$ if and only if 
$j:=max\{l;l\geq 2$ and $k \equiv 0$ mod $n_2 \cdots n_{l-1} \}\leq g$ and $m\geq \frac{k}{n_2 \cdots n_{j-1}}\overline{\beta}_j.$ Now recall
that $\overline{\beta}_{i+1}>n_i\overline{\beta}_{i}$ for $1\leq i \leq g-1$ and that $\overline{\beta}_{2}-n_1\overline{\beta_1}\geq e_2.$ This implies that for $3\leq j \leq g,$ we have  $\overline{\beta}_j- n_1 \cdots n_{j-1}\overline{\beta}_1>n_2 \cdots n_{j-1}(\overline{\beta}_{2}-n_1\overline{\beta}_1)\geq n_2 \cdots n_{j-1}e_2\geq e_1.$ So if 
$j\geq 3$ and $\kappa$ is a positive integer such that $m\geq \kappa\overline{\beta}_{j},$
we have $\frac{m-e_1}{n_1\beta_1}>\kappa n_2 \cdots n_{j-1},$ hence $q=[\frac{m-e_1}{n_1\beta_1}]\geq \kappa n_2 \cdots n_{j-1}.$ Therefore 
for $j\geq 3$, there are exactly $[\frac{m}{\overline{\beta}_{j}}]$ integers $\kappa \geq 1$ such that 
$m\geq \kappa \overline{\beta}_{j}$ and $\kappa n_2 \cdots n_{j-1}\leq q$, among them $[\frac{m}{n_j\overline{\beta}_{j}}]$  are $\equiv 0$ mod $(n_j).$\\
Similarly if $(q+1)n_1\beta_1+e_1\leq (q+1)\overline{\beta}_2,$ or equivalently $q\geq q_0,$ and if $\kappa$ is a positive integer 
such that $m\geq \kappa \overline{\beta}_2,$ we have $\kappa \leq \frac{m}{ \overline{\beta}_2}<q+1.$ Therefore if $q\geq q_0+1,$
we conclude that there are $\sum^g_{j=2}([\frac{m}{\overline{\beta}_j}]-[\frac{m}{n_j\overline{\beta}_j}])$
 empty $C_m^k$'s
 with $1\leq k \leq q.$
Moreover we have shown in the proof of theorem $4.9$ that $codim(C_m^0,\mathbb{C}^2_m)=codim(Cont^{>qn_1}(x_0)_m,\mathbb{C}^2_m)=
2+[\frac{m}{\beta_0}]+[\frac{m}{\beta_1}]$ if $m<(q+1)n_1\beta_1(resp.1+(q+1)(n_1+m_1)=1+[\frac{m}{\beta_0}]+[\frac{m}{\beta_1}]$
for $m\geq (q+1)n_1\beta_1).$Also note that $q_0\overline{\beta}_2<q_0n_1\beta_1+e_1<(q_0+1)n_1\beta_1+e_1\leq (q_0+1) \overline{\beta}_2\leq n_2\overline{\beta}_2<\overline{\beta}_3 \cdots$. Therefore for $q_0n_1\beta_1+e_1\leq m <(q_0+1)n_1\beta_1+e_1,$ we have $[\frac{m}{\overline{\beta}_2}]=q_0,[\frac{m}{n_2\overline{\beta}_2}]=[\frac{m}{\overline{\beta}_3}]= \cdots =0,$ so $N(m)=1,$
i.e. $C_m^0$ is irreducible.\\
Finally, assume that $qn_1\beta_1+e_1 \leq m < (q+1)n_1\beta_1+e_1$ with $q\geq 1$ and $q\leq q_0.$ Since 
$q_0<n_2,$ for $1\leq k \leq q$ we have $k\not \equiv0$ mod$(n_2)$ and $m \geq qn_1\beta_1+e_1>q\overline{\beta}_2,$
hence for $1\leq k \leq q,C_m^k=\emptyset$ and $C_m^0=Cont^{qn_1}(x_0)_m$ is irreducible.(The case $q=q_0$ was already known).So 
for $n_1\beta_1\leq m<(q_0+1)n_1\beta_1+e_1$, $C_m^0$ is irreducible.( Recall that for $1\leq m<q_0n_1\beta_1+e_1,$ the irreducibility of $C_m^0$ is already known).It only remains to check the codimensions of $C_m^0$ for $1 \leq m\leq q_0n_1\beta_1+e_1$. Here again we have seen in the proof of Theorem $4.9$ that $codim(C_m^0,\mathbb{C}^2_m)=codim(Cont^{>qn_1}(x_0)_m,\mathbb{C}^2_m)=:\delta^q(m)$ for any $q\geq 1$ and $qn_1\beta_1+e_1 \leq m <(q+1)n_1\beta_1+e_1$ and that 
$\delta^q(m)=$
$$2+[\frac{m}{\beta_0}]+[\frac{m}{\beta_1}] ~~for~~any ~~ q \geq 1 ~~and ~~ qn_1\beta_1+e_1 \leq m <(q+1)n_1\beta_1~~ $$
$$(q+1)(n_1+m_1)+1=1+[\frac{m}{\beta_0}]+[\frac{m}{\beta_1}] ~~for ~~q< q_0~~and ~~ (q+1)n_1\beta_1 \leq m <(q+1)\overline{\beta}_2~~ $$ $$(q+1)(n_1+m_1)+2=2+[\frac{m}{\beta_0}]+[\frac{m}{\beta_1}] ~~for ~~q< q_0~~and ~~ (q+1)\overline{\beta}_2\leq m <(q+1)n_1\beta_1 +e_1.$$
This completes the proof.

\end{dem}
In \cite{I}, Igusa has shown that the log-canonical threshold of the pair $((\mathbb{C}^2,0),(C,0))$ is 
$\frac{1}{\beta_0} +\frac{1}{\beta_1}.$ Here $(\mathbb{C}^2,0)$(resp.$(C,0))$) is the formal neighberhood of $\mathbb{C}^2$
(resp. $C$) at $0.$ Corollary .$4.10$ allows to recover corollary B of \cite{ELM} in this special case.

 \begin{cor} If the plane curve $C$ has a branch at $0$, with multiplicity $\beta_0$, and first Puiseux exponent $\beta_1,$
then 
$$ min_m \frac{codim(C_m^0,\mathbb{C}^2_m)}{m+1}=\frac{1}{\beta_0} +\frac{1}{\beta_1}.$$
\end{cor}
\begin{dem}
 For any $m,p\not=0$ in $\mathbb{N},$ we have $m-p[\frac{m}{p}]\leq p-1$ and $m-p[\frac{m}{p}]= p-1$ if and only if $m+1 \equiv 0$
mod $(p);$ so for any $m \in \mathbb{N},2+[\frac{m}{\beta_0}]+[\frac{m}{\beta_1}]\geq(m+1)(\frac{1}{\beta_0} +\frac{1}{\beta_1})$
and we have equality if and only if $m+1 \equiv 0$ mod $(\beta_0)$ and mod $(\beta_1)$ or equivalently 
$m+1 \equiv 0$ mod $(n_1\beta_1)$ since $n_1\beta_1$ is the least common multiple of $\beta_0$ and $\beta_1$.If not we have 
$1+[\frac{m}{\beta_0}]+[\frac{m}{\beta_1}]\geq (m+1)(\frac{1}{\beta_0} +\frac{1}{\beta_1}).$ Now if
$(q+1)n_1\beta_1 \leq m <(q+1)n_1\beta_1+e_1$ with $q \in \mathbb{N},$we have $(q+1)n_1\beta_1<m+1 \leq (q+1)n_1\beta_1+e_1<(q+2)n_1\beta_1,$
so $m+1 \not \equiv 0$ mod $(n_1\beta_1)$. If $(q+1)n_1\beta_1 \leq m <(q+1)\overline{\beta}_2$ with $q \in \mathbb{N}$ and $q<q_0$, then $(q+1)n_1\beta_1<m+1\leq (q+1)n_1\beta_1+e_1<(q+2)n_1\beta_1,$ so $m+1 \not \equiv 0$ mod $(n_1\beta_1).$ So in both cases, we have 
$1+[\frac{m}{\beta_0}]+[\frac{m}{\beta_1}]\geq (m+1)(\frac{1}{\beta_0} +\frac{1}{\beta_1}).$ The claim follows from corollary $4.10.$

 \end{dem}
It also follows immediately from corollary $4.10$

\begin{cor}Let $q_0 \in \mathbb{N}$ as in corollary $4.10.$ There exists $n_1\beta_1$ linear functions,
$L_{0}, \cdots ,L_{n_1\beta_1-1}$ such that $dim (C_m^0)=L_i(m)$ for any $m \equiv i$ mod $(n_1\beta_1)$ such that $m \geq q_0n_1\beta_1+e_1.$
\end{cor}
The canonical projections $\pi_{m+1,m}:C_{m+1}^0\longrightarrow C_m^0,m\geq 1,$ induce infinite 
inverse systems $$\cdots B_{m+1}\longrightarrow B_m \cdots\longrightarrow B_1$$
$$\cdots C_{(m+1)\kappa I}\longrightarrow C_{m\kappa I}\cdots\longrightarrow C_{(\kappa\beta_0\beta_1+e_1)\kappa I}\longrightarrow B_{\kappa\beta_0\beta_1+e_1-1}$$
and finite inverse systems $$C^j_{(\kappa\overline{\beta}_j-1)\kappa v}\longrightarrow C^j_{m\kappa v}\cdots\longrightarrow C^j_{(\kappa n_1 \cdots n_{j-1}\beta_1+e_1)\kappa v}\longrightarrow B_{\kappa n_1 \cdots n_{j-1}\beta_1+e_1-1}$$
for $2\leq j\leq g,$ and $\kappa \not \equiv 0$ mod $(n_j).$\\
We get a tree $T_{C,0}$ by representing each irreducible component of $C_m^0,m \geq 1,$ by a vertex $v_{i,m},1\leq i\leq N(m),$
and by joining the vertices $v_{i_1,m+1}$ and  $v_{i_0,m}$ if $\pi_{m+1,m}$ induces one of the above maps between the corresponding irreducible components. We represent the tree for the branch defined by $f(x,y)=(y^2-x^3)^2-4x^6y-x^9=0,$
 whose semigroup is $(4,6,15).$
 
\includegraphics[width=160mm,height=195mm]{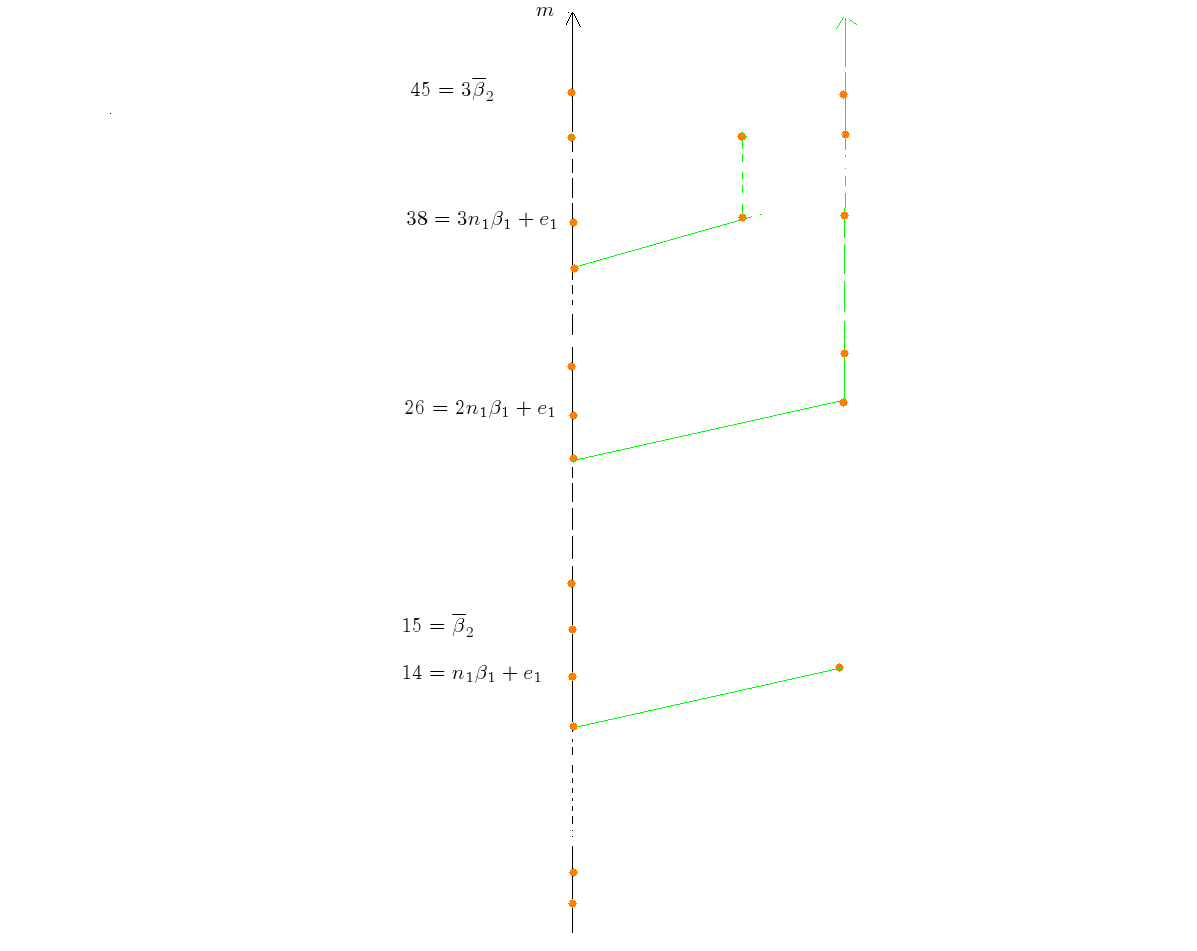}

This tree only depends on the semigroup $\Gamma.$\\
Conversely , we recover $\overline{\beta}_0, \cdots ,\overline{\beta}_g$ from this tree and $max\{m,codim(B_m,\mathbb{C}^2_m)=2\}= 
\overline{\beta}_0-1.$ Indeed the number of edges joining two vertices from which an infinite branch of the tree starts is 
$\beta_0\beta_1.$ We thus recover $\overline{\beta}_1$ and $e_1.$ We recover $\overline{\beta}_2-n_1\overline{\beta}_1, \cdots ,
\overline{\beta}_j-n_1 \cdots n_{j-1}\overline{\beta}_1, \cdots ,\overline{\beta}_g-n_1 \cdots n_{g-1}\overline{\beta}_1 ,$
hence $\overline{\beta}_2, \cdots ,\overline{\beta}_g$ from the number of edges in the finite branches.
\begin{cor}
Let $C$ be a plane branch with $g\geq 1$ Puiseux exponents. The tree  $T_{C,0}$ described above and $max\{m,dim~C_m^0=2m\}$
determine the sequence $\overline{\beta}_0, \cdots ,\overline{\beta}_g$ and conversely.
\end{cor}

\bibliographystyle{plain}

\thanks{Laboratoire de Math\'ematiques de Versailles,
 Universit\'e de Versailles-St-Quentin-en-Yvelines, 45 avenue des \'Etats-Unis,
78035 Versailles Cedex, France.}\\
\thanks{ Email address: mourtada@math.uvsq.fr}
 
\end{document}